\documentclass[12pt,leqno]{article}
\usepackage{amsmath, amsthm, amsfonts}
\begin{document}
\newtheorem{theo}{Theorem}
\newtheorem{prop}{Proposition}
\newtheorem{defi}{Definition}
\newcommand{\Hi}[1]{H_\nabla(#1)}
\def\P1{{\mathbb P}^1}
\def\C{{\mathbb C}}
\def\OP{\Omega_{\P1}^1(D)}
\def\O{{\cal O}}
\def\pr{{\bf Proof: }}
\begin{center}
{\LARGE\bf Meromorphic connections on $\P1$ and the multiplicity
of Abelian integrals\\}
\vspace{.25in}
{\large {\sc Hossein Movasati}
\footnote{Supported by MPIM-Germany
\\
\\
Math. classification: 14F05, 14F43}}
\end{center}
\begin{abstract}
In this paper we introduce the concept of Abelian integrals in differential
equations for an arbitrary vector bundle on $\P1$ with a meromorphic connection.
In this general context we give an upper bound for the numbers we are 
looking for.
\end{abstract}
Let $V$ be a locally free sheaf (vector bundle) of rank $\alpha$ on $\P1$ 
and $D=\sum_{i=1}^r m_ic_i$ be a positive divisor in $\P1$,  i.e. all $c_i$'s
are positive. We denote by $C$ the set of $c_i$'s. 
A meromorphic  connection $\nabla$ on $V$ with the pole divisor $D$ is 
a $\C$-linear homomorphism of sheaves 
$$
\nabla: V\rightarrow \OP\otimes_{\O_{\P1}}V
$$
satisfying the Leibniz identity
$$
\nabla(f\omega)=df\otimes \omega +f\nabla\omega,\  f\in \O_{\P1},\omega\in V
$$
 where $\OP$ is the sheaf of meromorphic 1-forms in $\P1$ with poles on $D$ 
(the pole order of a section of $\OP$ at $c_i$ is less than $m_i$).
For any two meromorphic connection  $\nabla_1$ and $\nabla_2$ with the same
pole divisor $D$, $\nabla_1-\nabla_2$ is a $\O_{\P1}$-linear map.

Let $t$ be the affine coordinate of $\C=\P1-\{\infty\}$, where 
$\infty$ is the point at infinity in $\P1$.
By Leibniz rule and by composing $\nabla$ with the holomorphic vector field 
$\frac{\partial}{\partial t}$ we can define:
$$
\nabla_{\frac{\partial}{\partial t}}: H^0(\P1, V(*\infty))\rightarrow H^0(\P1, V(D+*\infty))
$$
where $*\infty$ means that the
pole order at $\infty$ is arbitrary.
Since $\frac{\partial}{\partial t}$ is a holomorphic vector field in $\P1$ 
with a zero of multiplicity two at $\infty$, if 
$\omega\in H^0(\P1, V(*\infty))$ has a pole (resp. zero) of order $m$ at 
$\infty$ then $\nabla_{\frac{\partial}{\partial t}}\omega$ has a pole (resp. zero) of order $m-1$ 
(resp. $max\{2,m+1\}$) at $\infty$. If there is no confusion we write
$\nabla=\nabla_{\frac{\partial}{\partial t}}$. 

For any point $b\in \P1\backslash C$ we can find a frame $\{e_1,e_2,\ldots,
e_\alpha\}$  of holomrphic sections of $V$ in a neighborhood of $b$ such 
that $\nabla e_i=0\ \forall i$ and any other solution of
$\nabla \omega=0$ is a linear combination of $e_i$'s. Analytic continuations 
of this frame in $\P1\backslash C$ define the monodromy operator
\[
T:\pi_1(\P1\backslash C,b) \rightarrow GL(V_b)
\]
We say that $\nabla$ is irreducible if the action of monodromy on a non-zero
element of $V_b$ generates the whole $V_b$.
\\
Let $V^*$ be the dual vector bundle of $V$. There is defined a natural dual
connection $\nabla^*:V^*\rightarrow \OP\otimes_{\O_{\P1}}V^*$ on $V^*$ 
as follows
\[
<\nabla^*\delta,\omega>=d<\delta,\omega>- <\delta,\nabla\omega>, \ \delta\in V^*
\omega\in V 
\]
If $\{e_1,e_2,\ldots,e_\alpha\}$ is a base of flat sections in a neighborhood
of $b$ then we can define the dual of it as follows: $<\delta_i, e_j>=0$ if $i\not =j$ and $=1$ if $i=j$. We can easily check that $\delta_i$'s are flat
sections. The  associated  monodromy for $\nabla^*$ with respect to this basis
is just $T^*$, where $T^*$ is the composition of $T$ with the transpose operator.
We can also define a natural connection on 
$\wedge^k V=\{\omega_1\wedge\omega_2\wedge\cdots\wedge\omega_k
\mid \omega_i\in V\}$ with the pole divisor $D$ as follows:
\[
\nabla(\omega_1\wedge\omega_2\wedge\cdots\wedge\omega_k)=
\sum_{i=1}^k  \omega_1\wedge\omega_2\wedge\cdots
\widehat{\omega_i,\nabla\omega_i}\cdots \wedge\omega_k
\]
where $\widehat{\omega_i,\nabla\omega_i}$ means that we replace $\omega_i$ by
$\nabla\omega_i$.
\begin{prop}
\label{3nov01}
If the connection $\nabla$ over $V$ is irreducible then for any global 
meromorphic non zero section of $V$ with poles at $C\cup \{\infty\}$, 
say $\omega$, we have 
\begin{enumerate}
\item
$\{\nabla^i\omega\mid i=0,1,2,\ldots\}$ generates each fiber 
$V_b, b\in\P1\backslash C\cup \{\infty\}$;
\item
$\{\nabla^i\omega\mid i=0,1,2,\ldots,\alpha-1\}$ generates a generic 
fiber $V_b$. 
\end{enumerate}
\end{prop}
\begin{proof}
 If there exists a $\omega\in  H^0(\P1, V(*\infty))$ such that 
$\{\nabla^i\omega\mid i=0,1,2,\ldots\}$ does not generate $V_b$ then
there is a $\delta_b\in V_b^*$ such that 
\[
<\delta_b,\nabla^i\omega>=0,\ i=0,1,2,\ldots
\]
Consider the flat section $\delta$ passing through $\delta_b$. 
Since 
$\frac{\partial^i<\delta,\omega>}{\partial^i t}\mid_b=
<\delta_b,\nabla^i\omega>=0$,
we conclude that $<\delta,\omega>$ is identically zero. Since 
$\nabla$ is irreducible, we conclude that $\omega$ is the zero section.
\\
Now let us prove the second part. 
Let k be the smallest number such that for all non-zero 
$\omega\in  H^0(\P1, V(*\infty))$  $A=\omega\wedge\nabla\omega\wedge\cdots\wedge\nabla^k\omega$ is not 
identically zero. We want to prove that $k=\alpha-1$. 
 Fix a non zero $\omega$ with the property $A\wedge \nabla^{k+1}\omega=0$. 
Let $B=\P1-C\cup 
\{\infty\}\cup zero(A)$ and $V'$ be the vector 
bundle over $B$ generated by
$\omega,\nabla\omega,\cdots,\nabla^k\omega$. 
Since $A\wedge \nabla^{k+1}\omega=0$, $\nabla$ induces on
$V'$ a well-defined holomorphic connection. But this means that $V'_b$ is
invariant under monodromy. $\nabla$ is irreducible and so 
$k+1=dim(V'_b)=\alpha$.
\end{proof}

Every line bundle $L$ in $\P1$ is of the form $L_{a\infty}$, where $a$ is an 
integer and $L_{a\infty}$ is the line bundle associated to the divisor 
$a\infty$. We define $c(L)=a$ (Chern class).
According to Grothendieck decomposition 
theorem, every vector bundle $V$ on $\P1$ can be written as
$V=\oplus_{i=1}^\alpha L_i$, where $L_i$'s are line bundles. 
We define $c(V)=\sum_{i=1}^{\alpha}c(L_i)$. In view of Proposition ~\ref{3nov01}
the following definition is natural. 
\begin{defi}\rm
For any meromorphic global 
section of $V$ with poles at $C\cup\{\infty\}$ define its degree 
to be the sum of its pole orders. 
For any natural number $n$ let $\Hi{n}$ be the smallest number such that for 
all $\omega$ of degree $n$ the set
$\{\nabla^i\omega\mid i=0,1,2,\ldots, \Hi{n}-1\}$ generates each fiber 
$V_b, b\in\P1\backslash C\cup \{\infty\}$.
\end{defi}
 Of course
we have
\[
\Hi{n}\geq \alpha
\]
Let $V$ be a line bundle. In this case $\Hi{n}$ is the maximum multiplicity 
of a zero of a $\omega$ of degree $n$ minus one and so $\Hi{n}= n+c(V)+1$.
In general case we can only give an upper bound for $\Hi{n}$.

\begin{prop}
\label{26jan02}
Let $\nabla$ be an irreducible connection then
\[
\Hi{n}\leq (\alpha-1)(\sum m_i)+\alpha(n+1)-\frac{\alpha (\alpha-1)}{2}+ c(V)
\]
\end{prop}  
\begin{proof}
 For any global meromorphic section $\omega$ of $V$  with poles at $C\cup\{\infty\}$ we define 
$A=\omega\wedge\nabla\omega\wedge\cdots\wedge\nabla^{\alpha-1}\omega$. 
In Proposition ~\ref{3nov01} we proved that $A$ is
a nonzero global meromorphic section of $\wedge^\alpha V$. 
Let $m$ be the order of the pole of $\omega$ at $\infty$. The sum of pole orders
of $A$ at $C$ is at most $(\alpha-1)(\sum m_i)+\alpha (n-m)$. 
Each  $\nabla^i\alpha$ has a pole (resp. zero if $m-i$ is positive) of order $m-i$ 
at $\infty$ and so $A$ has order $m+m-1+m-2+\cdots+m-(\alpha-1)=m\alpha-\alpha(\alpha-1)/2$ 
at infinity. We conclude that the multiplicity of a zero of $A$ in $\P1\backslash C\cup\{\infty\}$ is less than
$$
(\alpha-1)(\sum m_i)+\alpha (n-m)+m\alpha-\alpha(\alpha-1)/2+c(V)
$$
If $b\in\P1\backslash C\cup\{\infty\}$ be a point such that
$\omega,\nabla
\omega\cdots,\nabla^i\omega, i\geq \alpha-1$ do not generate $V_b$ then $A$ has a zero of multiplicity $i-(\alpha-2)$. The proposition is proved.
\end{proof}
It does not seem to the author this upper bound to be the best one.
More precisely for any vector
bundle $V$ and divisor $D$ on $\P1$ can we find a meromorphic connection
$\nabla$ on $V$ with pole divisor $D$ such that $\Hi{n}$ is the above number?

Let $\delta$ be a flat section of $V^*$ in a small open set $U$ around $b$ 
and $\omega$ be a global meromorphic section of $V$ with poles at
$C\cup\{\infty\}$. From now on we use the notation $\int_{\delta}\omega$ 
instead of $<\delta,\omega>$.
Let us fix the number $n$ and  suppose that the degree of $\omega$ 
 is less than $n$. What is the maximum multiplicity of $\int_\delta\omega$  
at $t\in U$, say $\Hi{n,\delta}$?
Let $S(n)$ be  the vector space of meromorphic sections of 
$V$ with poles at $C\cup\{\infty\}$ and degree less than $n$. Since $S(n)$ is
a finite dimensional vector space, $\Hi{n,\delta_t}$ is  a finite number. 
\begin{prop}
\label{25nov01}
If $\nabla$ is irreducible then for all $t\in U$ 
\[
\Hi{n,\delta}\geq dim_\C S(n)-1
\]
The equality happens except for  a finite number of points in $U$.
\end{prop}
\begin{proof}
Let $\omega_1,\omega_2,\ldots,\omega_b$ be a basis for the vector space $S(n)$.
Consider  the determinant
\[
W_b(t)=det[\frac{\partial^i \int_\delta \omega_j}{\partial t^i}]_{b\times b}
\]
It is enough to prove that $W_b(t)$ is not identically zero. Let $a\leq b$
 be the smallest number such that $W_a(t)$ is identically zero. There exist holomorphic functions
 $p_i,i=1,2,\ldots,a$ in $U$ such that 
\[
A_a=\sum_{i=1}^{a-1} A_ip_i/p_a=0
\]
where $A_i$ is the $i$-th column of 
$[\frac{\partial^i \int_\delta \omega_j}{\partial t^i}]_{a\times a}$. This is $a$ equalities. If we act $\frac{\partial }{\partial t}$ to the $i$-th equation
and subtract the $(i+1)$-th equation we conclude that 
$[\frac{\partial^i \int_\delta \omega_j}{\partial t^i}]_{b-1\times b-1}.
[\frac{\partial (p_i/p_b)}{\partial t}]_{b-1\times 1}=0$. By hypothesis,
this implies that $[\frac{\partial (p_i/p_b)}{\partial t}]_{b_1\times 1}
\equiv 0$ or equivalently $\sum_{i=1}^{a}c_i\int_{\delta}\omega_i=0$, where 
$c_i$'s are constant. Since $\nabla^*$ is irreducible,  we have 
$\sum_{i=1}^{a}c_i\omega_i=0$ which is a contradiction.
\end{proof}
\begin{prop}
If $\nabla$ is irreducible then
\[
\Hi{n}=supp\{\Hi{n,\delta}\} 
\]
where $\delta$ runs through all flat sections of $V^*$ in $\P1\backslash C\cup\{\infty\}$.
\end{prop}
\begin{proof}
The proof is essentially stated in
Proposition ~\ref{3nov01}. 
If there exists a degree $n$ section $\omega$ such that 
$\{\nabla^i\omega\mid i=0,1,2,\ldots, p-1\}$ does not generate $V_b$ then
there is a $\delta_b\in V_b^*$ such that 
\[
\int_{\delta_b}\nabla^i\omega=0,\ i=0,1,2,\ldots,p-1
\]
Consider the flat section $\delta$ passing through $\delta_b$. 
We conclude that $\int_{\delta}\omega$ has multiplicity $p$ at $b$ and
hence $\Hi{n,\delta_b}\leq \Hi{n}$. The proof of the other part is similar.
\end{proof}
{\bf Regular Connections and Linear Equations:} 
Consider the connection $\nabla^*$ on $V^*$ as before and fix a 
trivialization 
map for $V^*$ around a singular point $c_i$. 
$\nabla^*$ is called regular  at $c_i$ if each flat section of $V^*$ in 
a sector 
with the vertex $c_i$ has at most a polynomial growth near $c_i$  
(see \cite{kul} p. 36 or 
\cite{ano} p. 8).
$\nabla$ is called regular if it is regular in all $c_i$'s.

Let $\omega_1,\omega_2,\ldots,\omega_\alpha$ be global meromorphic sections
of $V$ with poles at $C\cup\{\infty\}$. Let also $\delta_1,\ldots,\delta_\alpha$
be a base of flat sections of $V^*$ in a neighborhood of $b$.  The Wornskian function is defined as
follows
\[
W(t)=W(\omega_1,\cdots,\omega_\alpha)(t)=det[\int_{\delta_j}\omega_i]_{\alpha
\times\alpha}
\]
The division of two such functions is a one valued meromorphic 
function in $\P1\backslash C$ and by regularity of $\nabla^*$ we conclude 
that it extends meromorphically to the whole $\P1$. Fix an $\omega$.
By a similar argument as stated in Proposition ~\ref{25nov01} and by irreducibility of
$\nabla$ we know that
 $W(\omega,\nabla\omega,\ldots,\nabla^{\alpha-1}\omega)$ is not identically 
zero.  
The set $\{\int_{\delta}\omega\mid \delta \hbox{ is a flat section of } V^* \}$ is a base for the space of
solutions of the following linear equation
\begin{equation}
\psi: \left| \begin{array}{ll}
Y & \int_{\delta}\omega \\
Y' &  \int_{\delta}\nabla \omega \\
\vdots & \vdots \\
Y^{(\alpha)} &  \int_{\delta}\nabla^{\alpha} \omega
\end{array} \right|
\end{equation}
writing in other form
\begin{equation}
\label{bieshg}
\psi: Y^{(\alpha)}+\sum_{i=1}^{\alpha}(-1)^iP_i
Y^{(i)}=0
\end{equation}
where
$$
P_i=\frac{W(\omega,\nabla\omega,\ldots,
\widehat{\nabla^{\alpha-i}\omega},
\ldots,\nabla^{\alpha}\omega)}
{W},\ W=W(\omega,\nabla\omega,\ldots,\nabla^{\alpha-1}
\omega)
$$
Since $\int_{\delta_i}\tilde{\omega}$ has polynomial growth at the points of $C$, $\psi$ is regular therefore it must be Fuchsian i.e. $P_i$ has poles of
order at most $i$(see \cite{ano}). The union of poles of $P_i$'s is the singular
set of the Picard-Fuchs equation $\psi$. It has three type of singularities:
\begin{enumerate}
\item
$C$; in a $c_i\in C$ the solutions of ~(\ref{bieshg}) branch.
\item
$Z$ the zeros of $W$; In these singularities like regular points we have
a space of solutions of dimension $\alpha$. Note that
$P_1=\frac{\frac{\partial W}{\partial
t }}{W}$ and so neither of these points is regular. For this reason in $\cite{ano}$ these are called apparent singularities.
For a zero $b$ of $W$
we can find a flat section $\delta$ of $V^*$ 
such that $\int_{\delta}\omega$ has multiplicity greater than $\alpha$ at
$b$.
\item
$\infty$; 
Let $m$ be the order of the pole of $\omega$ at $\infty$. 
 The solutions of
~(\ref{bieshg}) in a neighborhood of $\infty$ are meromorphic functions with poles of order  at most $m$ at $\infty$.
 \end{enumerate}
Since $P_1=\frac{\frac{\partial W}{\partial t }}{W}$ we have
$$
Res(P_1dt, t=c)=mul(W, t=c),\ c\in C\cup Z
$$
Now Consider a  regular linear equation $\psi$ with singularities at 
$C\cup Z\cup \{\infty\}$ and suppose that it has apparent singularities in $Z$ and a singularity of type 3 in $\infty$. Furthermore assume that $\psi$
 has the same
  monodromy representation like as of $\nabla^*$.
\begin{prop}
$\psi$ is obtained 
by a meromorphic global section of $V$ with poles at $C\cup\{\infty\}$.
\end{prop} 
\begin{proof}
Consider in a neighborhood of $b$ a base of flat sections $\delta_i$ of $V^*$
and a base $e_i,i=1,2,\ldots, \alpha$ for solutions of $\psi$ such that the monodromy representation of
the both $\nabla^*$ and $\psi$ with respect to these bases is the same
Define a section of $V=(V^*)^*$ as follows
\[
\omega (\delta_t)=e_i(t) 
\]
this is a one valued holomorphic section of $V$ in $\P1\backslash C\cup\{\infty\}$. Since $\psi$ and $\nabla^*$ are regular, $\omega$ 
extends meromorphically to $C$. 
\end{proof}

Let $\phi_b$ be the maximum multiplicity of solutions of $\psi$ at $b$. If $b$
is a regular point of $\psi$ then $\phi_b=\alpha-1$ and if it is an apparent
singularity of $\psi$ then $\phi_b\geq \alpha$. In the last case by 
definition of $W$ we can see that $W$ has a zero of order at least $\phi_b-
(\alpha-1)$ at $b$ and by $P_1=\frac{\frac{\partial W}{\partial
t }}{W}$ we have
\[
\phi_b\leq Res(P_1dt ,t=b)+(\alpha-1)
\]
{\bf Remark:} Let us choose a trivialization of $V$ in a small
disk $D$ around  a singular point $c_i$ of the connection $\nabla$, $V\mid_{D}\cong D\times \C^\alpha$, and a coordinate $z$ in
$D$. In this coordinate we can 
write $\nabla v=\frac{\partial v}{\partial z} +
\sum_{j=1}^{m_i}\frac{C_j}{z^j}v+ A(z)v$, where $v$
is a holomorphic vector in $D$, $C_j, 1\leq j\leq m_i$ (resp. $A(z)$) is a constant (resp.
 holomorphic in $z$) matrix. $C_1$ is called the residue of the connection at
$c_i$. Now we can apply the Levelt's theory (see
\cite{ano} Section 1, 2.2) to understand the local theory of this 
connection.  
      \newcommand{\dis}[2]{R^{#1}f_*#2}
       \def\RR{{\cal R}}
       \newcommand{\form}[1]{\Omega^1(#1D)}
       \newcommand{\cfv}[1]{\{\delta_t\}_{t\in #1}}

{\bf Lefschetz Pencil:}
Let $M$ be a projective compact complex manifold of dimension
two, $\{M_t\}_{t\in\P1}$ a pencil of hyperplane sections of $M$
and $f$ the meromorphic function on $M$ whose level sets are
$M_t$'s (see \cite{lam}). We set $\RR$ the indeterminacy points of $f$,
$L_t=M_t-\RR$, $C=\{c_1, c_2,c_3,\ldots,c_r\}$ the set of critical
values of $f$, $\beta=dim(H^1(L_t,\C))$ for a $t\in\P1-C$ and $\C[t]$ the 
ring of polynomials in $t$. Since $f\mid_{M-\RR}$ is a $C^\infty$ fibration over $\P1-C$ (see
\cite{lam}), $\beta$ is independent of $t$. We assume that
\begin{enumerate}
\item
The axis of the pencil
intersects $M$ transversally. This is equivalent to this
fact that in a coordinate system $(x,y)$ around each indeterminacy point of $f$ we can write $f=\frac{x}{y}$;
\item
The critical points of $f$ are isolated;
\item
The pole divisor $D=M_\infty$ of $f$ is a regular fiber, 
i.e. $\infty\not\in C$. 
\end{enumerate}

We define $\Omega^i(*D)$ to be the set of meromorphic $i$-forms in $M$ with poles of arbitrary order along $D$. 
The set $\tilde{H}=\cup_{t\in B}H^1(L_t,\C)$, where $B=\P1-C$,
has a natural structure of a complex manifold and the natural
projection $\tilde{H}^1\rightarrow B$ is a holomorphic vector
bundle which is called the cohomology vector bundle.
The sheaf of holomorphic sections of $\tilde{H}$ is also denoted by $\tilde{H}$.
In what follows when we consider $f$ as a holomorphic function we mean
its restriction to $M-\RR$. Let $\C_{M-\RR}$ be the sheaf of constant functions in $M-\RR$ and  $\dis{1}{\C_{M-\RR}}$ be the $1$-th direct image of the
sheaf $\C_{M-\RR}$ (see \cite{gra2}). Any element of
$\dis{1}{\C_{M-\RR}}(U)$, $U$ being an open set in $B$, is a
holomorphic section of the cohomology fiber bundle map. It is easy to verify that
\[
\tilde{H}\cong \dis{1}{\C_{M-\RR}}\otimes_{\C}\O_{\P1},
\hbox{ in } B
\]
Now let us introduce the Gauss-Manin Connection on
$\tilde{H}$. Consider a holomorphic coordinate $(t,0)$ in $U$,
a small open disk in $\P1$. The Gauss-Manin connection is defined
as follows:
\[
\nabla:
\tilde{H}(U)\rightarrow \Omega^1_{\P1}\otimes_{\O_U}
\tilde{H}(U)
\]
\[
\nabla(g\otimes c)=dg \otimes c,\ c\in \dis{1}{\C_{M-\RR}}(U),\ 
g\in \O_{\P1}(U)
\]
The sheaf of flat sections of $\nabla$ is $\dis{1}{\C_{M-\RR}}$.
Let $\frac{\partial}{\partial t}$ be a vector field in $U$. 
We write $\nabla_{\frac{\partial}{\partial t}}=
\frac{\partial}{\partial t}\circ\nabla$. ($\nabla_
{\frac{\partial}{\partial t}}(g\otimes c)= 
\frac{\partial g}{\partial t}\otimes c$).

In the same way we can define the cohomology fiber bundle $\tilde{H}_c$ of 
compact fibers $M_t$. Since $\tilde{H}_c$ is
a $\nabla$-invariant sub-vector bundle of $\tilde{H}$, we have the restriction 
of $\nabla$ to $\tilde{H}_c$ which we denote it again by $\nabla$.

Let $\omega$ be a meromorphic 1-form in $M$ with poles along some fibers of $f$.
 Let also $\cfv{\P1-C},\delta_t\subset L_t$ be a a continuous 
family of cycles. The Abelian integral $\int_{\delta_t}\omega$ appears in
the deformation $df+\epsilon\omega$ of $df$ inside holomorphic foliations 
(differential equations) and it is related to the number of limit cycles which
born from the cycles $\delta_t$ (see \cite{hos1}). 
The pair $(\tilde{H},\nabla)$ is defined in $\P1-C$ and in order to be in the
context of this paper we may be interested to prove:
\begin{prop}
\label{13jan02}
Under the assumptions 1,2,3, there is a vector bundle
$V$, a sub-vector bundle $\bar{V}\subset V$ and a meromorphic connection on $V$ with poles in $C$ 
$$
\nabla: V\rightarrow \Omega^1_{\P1}(D)\otimes_{\O_{\P1}}V,\ D=\sum m_ic_i
$$
such that: 
\begin{enumerate}
\item 
$\bar{V}$ is $\nabla$-invariant.
\item
$(V,\nabla)$ (resp. $(\bar{V}, \nabla)$) coincides with $(\tilde{H}, \nabla)$( resp. $(\tilde{H}_c,\nabla)$) in $\P1-C$;
\item 
The Brieskorn lattices ( Petrov module in the
context of differential equations) of $f$ $H'=\frac{\Omega^1(*D)}{df\wedge\Omega^0(*D)+d\Omega^0(*D)}$ is $\C[t]$-isomorphism to the module
of global sections of $V$ with poles of arbitrary order at $\infty$.
\end{enumerate}
\end{prop}
This is a task which is done in details in \cite{hos}.   
If the singularities of $f$ are non-degenerated, i.e. in a holomorphic 
coordinate $(x,y)$ around a singularity $p_i$ we can write $f=f(p_i)+x^2+y^2$, 
then all the $m_i$'s are equal to one. In other words $\nabla$ is
logarithmic.

The pair $(V,\nabla)$ is not irreducible but if $H^1(M,\C)=0$ and 
$f$ satisfies 1, 2, 3 and has
non-degenerated singularities with distinct images then $(\bar{V}, \nabla)$ is
irreducible (see \cite{lam} 7.3). The following proposition justifies the
use of $(\bar{V},\nabla)$ instead of $(V,\nabla)$.
\begin{prop}
\label{16jan02}
For an $\omega$ a meromorphic 1-form in $M$ with poles of order at most $n$
along $D$, the integral $\int_{\delta_t}\omega$ is a polynomial
of degree $n$. $\nabla_{\frac{\partial}{\partial t}}^i\omega, i>n$ 
restricted to each fiber
has not residues in $\RR$ and hence is a meromorphic section of $\bar{V}$.
\end{prop}
\begin{proof}
We have $p(t):=\int_{\delta_t}\omega=
t^n\int_{\delta_t}\frac{\omega}{f^{n}}$. Since the 
1-form 
$\frac{\omega}{f^{n}}$ has not pole along $D$,
 $\frac{p(t)}{t^n}$ has finite growth at $t=\infty$. Since $p(t)$ is
holomorphic in $\C$ (even in the points of $C$), we conclude that
$p(t)$ is a polynomial of degree at most
$n$. The second part is a direct consequence of the first one and
the fact that $\frac{\partial}
{\partial t}\int_{\delta_t}\omega=\int_{\delta_t}\nabla_{\frac{\partial}{\partial t}}\omega$.
\end{proof}

I tried to study the maximum multiplicity of Abelian integrals in the context
of meromorphic connections. My motives were the paper \cite{mar} and also a
paper of mine, where the extension of cohomology vector bundles and their
connections to the critical values of a meromorphic function is discussed.
The upper bound obtained in Proposition ~\ref{26jan02} seems to be
far from the best one (at least for Guass-Manin connections). Some works 
in Differential equations (see \cite{hor}) suggest that the number
$H_{\nabla}(n)$ must be very sensitive with respect to $\nabla$. 

Here I would like to express my thanks to Max-Planck institute for
hospitality. Thanks go also to C. Hertling, Y. Holla and S. Archava for many useful conversations. 

\smallskip
\leftline{Hossein Movasati} 
\leftline{IMPA, 22460-320, Rio de Janeiro, Brazil} 
\leftline{IPM, P.O.Box : 19395-5746,Tehran-Iran} 
\leftline{E-mail: hossein@impa.br, movasati@ipm.ir}
 
\end{document}